# THE EXPECTED DURATION OF RANDOM SEQUENTIAL ADSORPTION

BY AIDAN SUDBURY

*Monash University*

When gas molecules bind to a surface they may do so in such a way that the adsorption of one molecule inhibits the arrival of others. We consider random sequential adsorption in which the empty sites of a graph are irreversibly occupied in random order by a variety of types of "particles." In a finite region the process terminates when no more particles can arrive. A universal asymptotic formula for the mean duration is given.

**1. Introduction.** In this paper we consider a process known as *random sequential adsorption* (*RSA*) in which "particles" are deposited on a graph successively, subject to a blocking effect, whereby the presence of a particle at a particular set of sites may prevent the arrival of other particles at the same or nearby sites. The principal motivation comes from the adsorption of gas molecules onto a surface. Initially the graph is empty and, at rate 1, independently of all other particles and sites, each type of particle attempts to occupy each suitable set of sites. No particles leave the graph, so that once the arrival of a particle is blocked it is blocked forever. Thus, the process terminates on any finite set when no further particles can arrive. The most common types of RSA which have been considered are those in which single particles (monomers) arrive at single sites, but inhibit the arrival of neighboring particles, or when pairs of particles (dimers) can only arrive where pairs of positions are free. There are fairly general treatments of these processes in Penrose and Sudbury [2], a review article by Evans [1] and more recently the series of surveys in Volume 165 of the journal *Colloids and Surfaces* (2000).

We wish to consider RSA processes on $Z^d$ which are sufficiently flexible to allow a variety of possible changes at any site. Thus, we do not restrict ourselves to monomers and dimers or to any particular ways in which the arrival of one type may prevent the arrival of another type. To emphasize that there is no restriction on the shape of a particle, the term "config" is used in this paper. A configuration $T_i$ is defined by a finite set of sites which must include the origin as "anchor." Any translation of this set is called a "config" of type $T_i$. A config arrives at a site $x$ when its anchor attempts to arrive at $x$. If no neighbors of the config have arrived beforehand, then the arrival will be successful. We prove the following:

THEOREM 1. *Suppose we have an RSA process in which $k$ different types of config attempt to arrive at each site, each independently at rate 1. Consider a subset $A$ of $Z^d$ consisting of $n$ sites. Call the expected time until the process terminates on the $n$ sites $E_{\text{term}}(A)$. Call the set of subsets of $Z^d$ of $n$ sites $S_n$. Define $p$ to be the limiting probability as $t \to \infty$ that an arrival at a site will be successful. If $p > 0$, then, as $n \to \infty$,*

$$\max_{A \in S_n} \left| 1 + \frac{1}{2} + \frac{1}{3} + \cdots + \frac{1}{kn} + \ln(p) - E_{\text{term}}(A) \right| \to 0.$$

Exact values of $p$ are not generally known, but in 1-d the dimer process has $p = e^{-2}$ and the annihilation process $p = e^{-1}$ (see, e.g., Sudbury [3]). In the case of dimers arriving on $Z$, it would be the limiting probability that, given that a dimer had not arrived at the pair of sites 0, 1, dimers already occupied pairs $-2$, $-1$ and 2, 3. Simulations of the dimer process on a set of 400 sites gave a mean that differed from the asymptotic estimate of 4.567 by approximately 0.02. With this number of sites it can be seen that the $\ln(p) = -2$ term plays an important role. It will be clear that the variance of the expected time is bounded by $\sum 1/n^2 = \pi^2/6$.

**2. The expected time for the process.** If there was no interference between configs, then the expected time until all had arrived would be $1 + 1/2 + \cdots + 1/(kn)$. If the probability of the successful arrival of a config was $p$ independent of all other configs, then the expected time until the last arrival would be $\sum_{r=0}^{kn-1} p(1-p)^r \{\sum_{i=r+1}^{kn} i^{-1}\}$ which, asymptotically, tends to $1 + 1/2 + \cdots + 1/(kn) + \ln(p)$, as will be shown at the end of this section. The chief idea of the proof given here is that a subset of sites of $O(n)$ do behave as if they are almost independent when $n$ is large. This is because, although the number of neighbors at a distance $l$ in $Z^d$ grows as $l^d$, the influence of one site on another falls off as $1/l!$. Further, the "bulk of the action" in the series $1 + 1/2 + \cdots + 1/(kn)$ occurs at the beginning, that is, for the later stages of the RSA process.

It should be noted that the bounds given in the following lemmas do not depend on the positions of a set of configs, only on their number.

Assume we have a finite set of types of finite configurations of sites $\{T_1, T_2, \ldots, T_k\}$. [E.g., in 2-d we might have the two configurations $T_1 =$ $\{(0,0),(0,1),(1,1)\}$, $T_2 = \{(0,0),(2,0)\}$.] The origin is chosen to be the anchor. Any translation of a member of the set of types is a "config" [e.g., the set of points $(m-1,n-1)$, $(m-1,n)$, $(m,n)$ would be a config of type $T_1$]. The anchor translates with the config. For each config $C$ there is a single arrival time $t(C)$ which is a mean 1 exponential random variable, which is independent of the arrival times of all other configs. We may thus construct the process via a Harris diagram in the usual way.

Two configs $C_1, C_2$ are said to be neighbors if the successful arrival of one of them must block the successful arrival of the other. If a config has not arrived and one of its neighbors has made a successful arrival, then it enters a blocking mode. In a finite region, once all arrivals have taken place or all unarrived configs are in a blocking mode, the process terminates.

Assume that initially all the sites of $Z^d$ are unoccupied so that the process is translationally invariant. Suppose the process on $Z^d$ has been running for time $t$. Define $p_i(t)$ to be the probability that a particular config of type $T_i$ is not blocked at time $t$ given that it has not arrived. Since once any config is in blocking mode, it remains so, this probability is monotonic decreasing and $\lim_{t \to \infty} p_i(t) = p_i$.

In a realization, a config $C_0$ is said to "affect" another config $C$, if there is a sequence of configs $C_0, C_1, C_2, \ldots, C_l = C$ such that $t(C_i) < t(C_j)$, $i < j$, and $C_i$ is a neighbor of $C_{i+1}$. (Note that $C_0$ "affecting" $C$ does not imply that the arrival of $C_0$ either ensures the blocking of $C$ or its successful arrival. However, $C_0$ "affecting" $C$ is a necessary condition for the arrival or nonarrival of $C_0$ to affect $C$ in the sense that a realization with arrivals identical to the first, except for the nonarrival of $C_0$ before $C$, will differ from the first in the successful arrival or nonarrival of $C$.) All distances are $L_1$. Let $D_b$ be the largest distance from the anchor of one config to the anchor of a neighbor. An $m$-neighborhood of a config is the set of configs whose anchors are distance $\leq m$ from the config's own anchor. We put $N = nk$.

LEMMA 2.  *Define $A_m(C)$ to be the event that a config $C$ is affected by at least one config outside its $m$-neighborhood. Then $P(A_m(C)) < N^{-m_0}$ for any positive integer $m_0$ when $m = [D_b N^{1/2d}]$ and $N$ is large enough.*

PROOF.  The probability that in a config path of length $l$ the events are ordered in time is $= 1/l!$. Let $N_b$ be the maximum number of neighbors of any config, then the number of paths of length $l$ is $< (N_b)^l$. To be distance $m$ away, the path must be at least length $m/D_b$. Thus, the probability a config is affected by another at distance $m$ or greater is

$$P(A_m(C)) < \sum_{k=[m/D_b]}^{\infty} \frac{N_b^k}{k!} < 2\frac{N_b^{m/D_b}}{([m/D_b])!} < \left[\frac{(eN_b)^{2d}}{N}\right]^{N^{1/2d}/2d} < N^{-m_0}$$

for any positive integer $m_0$ when $m = [D_b N^{1/2d}]$ and $N$ is large enough. □

If there are no config paths from outside the $m$-neighborhood, it is clear that events outside the neighborhood cannot influence the config, thus, we have the following:

LEMMA 3. *Let $C$ be the event that the config $C$ arrives successfully. Let $E$ be an event defined by sites outside the $m$-neighborhood of $C$. It follows that $P\{C|E, A_m^c(C)\} = P\{C|A_m^c(C)\}$.*

The number of configs in the $m$-neighborhood is called $N_m$.

LEMMA 4. *In $d$-dimensions, $N_m < M\sqrt{N}$ where $m = [D_b N^{1/2d}]$ and $M$ is a fixed constant.*

PROOF. The set of sites with $L_1$ distance from the origin $\leq m$ is contained in the $d$-dimensional cube with each co-ordinate lying in $[-m, m]$. The number of sites is thus less than $(2m+1)^d$ and the number of anchors at each site is $k$. Putting $m = [D_b N^{1/2d}]$, the results follows. □

Let $Q_N$ be a set of $N$ configs. In any realization of the RSA process a subset of $r$ configs will be the last of the $N$ to arrive. Define the event

$$E(r, Q_N, m) = \{\text{Each of the last } r \text{ configs to arrive is}$$
$$\text{the last of its } m\text{-neighborhood}\}.$$

Note that arrival here does not imply successful arrival.

LEMMA 5. *If $N$ is large enough so that $rN_m < N$, then*

$$P\{E(r, Q_N, m)\} > 1 - \frac{r^2 N_m}{N} > 1 - M\frac{r^2}{\sqrt{N}}.$$

PROOF. Call the last $r$ configs to arrive $C_1, C_2, \ldots, C_r$. There are at least $N$ configs available to be $C_r$, and one of them must be the last of at most $N + rN_m$ configs. $C_r$ arrives after $C_1, C_2, \ldots, C_{r-1}$, and since each $C_i$ is also the last of its $m$-neighborhood, $C_r$ must arrive after the union of these $m$-neighborhoods. Thus, the probability of this occurring is $> N/(N + rN_m)$. Conditional on this event, we now calculate the lower bound of the probability for $C_{r-1}$ to be the last of its $m$-neighborhood. There are at least $N - N_m$ of the original $N$ available to be $C_{r-1}$ and one of them must be the last of at most $N + (r-1)N_m$ configs. Thus, the conditional probability that the second last is the last of its $m$-neighborhood is $> (N - N_m)/(N + (r-1)N_m)$. Continuing the argument, we obtain

$$P\{E(r, Q_N, m)\} > \frac{N(N - N_m) \cdots (N - (r-1)N_m)}{(N + rN_m)(N + (r-1)N_m) \cdots (N + N_m)}.$$

Now,

$$\frac{N - kN_m}{N + (r-k)N_m} > 1 - \frac{rN_m}{N} \Rightarrow P\{E(r, Q_N, m)\} > \left(1 - \frac{rN_m}{N}\right)^r,$$

from which the result follows using Lemma 4. □

Define $p = \sum_{j=1}^k p_j / k$.

LEMMA 6. *Given none of the last $r+1$ configs to arrive is in the $m$-neighborhood of any other, then the probability all except the first of them are blocked $= B^{r+1} p(1-p)^r$, where $e^{-N^{-m_0}} < B < e^{N^{-m_0}}$.*

PROOF. First, note that the assumptions imply that each of the last $r+1$ configs to arrive is the last of its $m$-neighborhood. Call the last $r+1$ configs $C_1, C_2, \ldots, C_{r+1}$ and define $E_i, i = 1, \ldots, r+1$, to be the event that $C_i$ is blocked. We now show that the dependence of each event on the others is small, because being affected by events outside the $m$-neighborhood is small:

$$P\left(E_i \Big| \bigcap E_j / E_i\right)$$
$$= P\left(E_i \Big| A_m(C_i), \bigcap_j E_j / E_i\right) P\left(A_m(C_i) \Big| \bigcap_j E_j / E_i\right)$$
$$+ P\left(E_i \Big| A_m^c(C_i), \bigcap_j E_j / E_i\right) P\left(A_m^c(C_i) \Big| \bigcap_j E_j / E_i\right)$$
$$= P\left(E_i \Big| A_m(C_i), \bigcap_j E_j / E_i\right) P(A_m(C_i)) + P(E_i | A_m^c(C_i)) P(A_m^c),$$

using Lemma 3 and the fact that the event $A_m(C_i)$ is determined by events in the $m$-neighborhood of $C_i$. The probability of $C_i$ being blocked unconditioned on the status of the other $r$ configs differs only in the first term above, which is less than $P(A_m)$. Thus, putting the conditional probability of being blocked as $\epsilon_i + (1 - p_i)$, where $\epsilon_i$ may be negative, and averaging over the $\{p_i\}$, we obtain $\epsilon + (1-p), |\epsilon| < P(A_m)$. This equals $B(1-p)$ where, using Lemma 2, we have $e^{-N^{-m_0}} < B < e^{N^{-m_0}}$. □

The values of $B$ will depend on time and on the positions of the configurations. In the interests of simplicity of presentation, the values of $B$ are not always the same in what follows, but they all satisfy the above inequality.

LEMMA 7. *The probability the last $r$ arrivals are blocked is $< \beta^r$ for fixed $\beta < 1$ when $r < kn/2N_m$. For $r > kn/2N_m$, this probability is $< \beta^{kn/2N_m}$.*

PROOF. Each site has $k$ configs anchored at it, and, thus, the number of sites represented in the last $r$ must be at least $r_k = [r/k]$. Now, select a subset of the $n$ sites of size at least $n/N_m$ consisting of sites none of which are in the $m$-neighborhood of each other. Assume $r_k < n/2N_m$. The number of sites from the subset in the last $r_k$ is a hypergeometric random variable which, however, stochastically dominates a binomial random variable. Even assuming that all of the previous $r_k - 1$ sites selected have been from the subset size $n/N_m$, the probability the $r_k$th is from the subset is $> (n/N_m - r/k)/n > 1/2N_m$ by assumption. Thus, the number of sites dominates a random variable $X \sim \text{Binom}(r_k, 1/2N_m)$. The conditional probability of any particular one of these sites being blocked is $p_{\text{block}} < 1 - p + P(A_m) < 1$, using arguments from the last lemma. Thus, the probability of all of these being blocked is less than

$$E\{p_{\text{block}}^X\} = [1 - (1 - p_{\text{block}})/2N_m]^{r_k} < \beta^r,$$

from which the first part of the result follows. Clearly, as $r$ increases the probability can only decrease. □

In the next lemma we shall write $E(r, N, m)$ instead of $E(r, Q_N, m)$ since the bound does not depend on the positions of the sites in $Q_N$.

LEMMA 8. *Assume that $n$ is large enough that $kn/2N_m = N/2N_m > n^{1/7}$. Then the probability the last change in the RSA process on a set of $n$ sites occurs at arrival $N - r$ is $< \beta^{n^{1/7}}$ for $r > n^{1/7}$. When $r < n^{1/7}$ the probability equals $B^{r+1}p(1-p)^r + \epsilon_r$, where*

$$|\epsilon_r| < M\frac{r(r+1)}{\sqrt{N}}$$

*and $e^{-N^{-m_0}} < B < e^{N^{-m_0}}$.*

PROOF. The first part follows directly from Lemma 7. Suppose now that $r < n^{1/7}$. If the event $E(r+1, N, m)$ occurs, then the last $r+1$ arrivals are at sites not within $m$ of each other. The conditional probability the last change is at $kn - r$ is $B^{r+1}p(1-p)^r$ by Lemma 6. We then use Lemma 5 to get the bound on $\epsilon_r$. □

We are now in a position to prove Theorem 1.

PROOF OF THEOREM 1. The times between arrivals are independent of the order of arrival of the configs so that, if exactly $j$ configs have arrived, the expected time until the next arrival is $1/(kn - j)$. The expected time until the $(kn - j)$th arrival is thus $\sum_{i=j+1}^{kn} 1/i$. It follows that the expected time until the last successful arrival is

$$E_{\text{term}} = \sum_{r=0}^{[n^{1/7}]} \{Bp(B(1-p))^r + \epsilon_r\}\left\{\sum_{i=r+1}^{kn} \frac{1}{i}\right\} + \sum_{r=[n^{1/7}]+1}^{kn-1} \eta_r \left\{\sum_{i=r+1}^{kn} \frac{1}{i}\right\},$$

where each $\eta_r < \beta^{n^{1/7}}$. As $n \to \infty$, the last term tends to 0. Using Lemma 8, the term in $\epsilon$ breaks up into

$$\sum_{r=0}^{[n^{1/7}]} \epsilon_r \left\{\sum_{r+1}^{[n^{1/7}]} \frac{1}{i}\right\} < M \sum_{r=0}^{[n^{1/7}]} \frac{1}{i}\sum_{0}^{i} \frac{r(r+1)}{\sqrt{N}} < \frac{M}{2\sqrt{N}} \sum_{r=0}^{[n^{1/7}]} i^2 \to 0$$

and

$$\sum_{r=0}^{[n^{1/7}]} \epsilon_r \left\{\sum_{[n^{1/7}]}^{kn} \frac{1}{i}\right\} < M \sum_{r=0}^{[n^{1/7}]} \frac{r(r+1)\log(kn)}{\sqrt{N}} \to 0.$$

Rearranging the last equation gives

$$\sum_{i=1}^{kn} \frac{1}{i} \sum_{r=0}^{\min(i-1,[n^{1/7}])} \{Bp(B(1-p))^r\} = Bp\sum_{i=1}^{kn} \frac{1}{i}\frac{1 - (B(1-p))^{\min(i,[n^{1/7}]+1)}}{(1 - B(1-p))}.$$

As $n \to \infty$, $B \to 1$, and this expression equals

$$\sum_{i=1}^{kn} \frac{1}{i} + \ln(1-(1-p)) - (1-p)^{[n^{1/7}]+2}\sum_{i=[n^{1/7}]+2}^{kn} \frac{1}{i} + \sum_{i=[n^{1/7}]+2}^{\infty} \frac{(1-p)^i}{i},$$

and the last terms tend to 0 as $n \to \infty$. □

**3. The annihilating process.** In RSA empty sites become filled as the process proceeds, whereas in annihilating processes, the sites begin occupied, and then some become unoccupied as time passes. In many cases a simple swapping of the roles of occupied and unoccupied turns one sort of process into the other. Annihilating processes have been considered in Sudbury [3] and Penrose and Sudbury [2]. Initially all sites of a graph are occupied. Between each pair of adjacent sites occur annihilation events at exponential rate 1. If one of the sites is unoccupied, nothing happens. If both are occupied, then one of the two sites, chosen at random, becomes unoccupied. (By swapping the roles of "occupied" and "unoccupied," we have an RSA process in the sense of the previous sections, where a pair 00 can be replaced by either 01 or 10.) If we start with $n$ particles in a line and they evolve as an annihilating process, then at some time $T_n$ the process will stop as all remaining particles are separated from each other. We solve for the distributions of the times in the sense that we find a generating function for them. Define

$$F_n(t) = P(T_n < t). \tag{1}$$

Then it is clear that $F_0(t) = F_1(t) = 1$. Further, if there are two separated blocks of $m$ and $n$ particles, then

$$F_m(t)F_n(t) = P(T_m, T_n < t).$$

In a block of $n$ particles, the end particles are annihilated at rate $1/2$, and the interior particles are annihilated at rate 1, so the time until the first break is exponential with rate $n-1$. Thus,

$$F_n(t) = \int_0^t \sum_{r=1}^{n-1} F_r(t-s)F_{n-1-r}(t-s)e^{-(n-1)s}\,ds, \qquad n \geq 2. \tag{2}$$

Put $u_n(t) = F_n(t)e^{nt}$, then

$$u_n(t) = e^t \int_0^t \sum_{r=1}^{n-1} u_r(t-s)u_{n-1-r}(t-s)\,ds, \qquad n \geq 2. \tag{3}$$

With $U(x,t) = \sum_{n=1}^\infty u_n(t)x^n$, we have

$$\sum_{n=2}^\infty \sum_{r=1}^{n-1} u_r u_{n-r-1} x^n = u_1 x + (u_1 u_1 + u_2)x^2 + \cdots = U + U^2.$$

Thus, multiplying equation (3) by $x^n$ and adding from 2 to $\infty$, we obtain

$$U - xe^t = xe^t \int_0^t U(t-s) + U^2(t-s)\,ds,$$

or

$$Ue^{-t} - x = x\int_0^t U(s) + U^2(s)\,ds,$$

and differentiating gives

$$e^{-t}\frac{dU}{dt} - Ue^{-t} = x(U + U^2).$$

Putting $v = 1/U$, we obtain

$$\frac{dv}{dt} + (xe^t + 1)v = -xe^t,$$

and using an integrating factor gives

$$v = -1 + \frac{1}{x}e^{-t} + c(x)e^{-(xe^t+t)},$$

where $c(x)$ is a function to be determined from the boundary conditions. $F_1(0) = 1$ and $F_i(0) = 0$, $i > 1$. Thus, $U(x,0) = x$ and

(4) $$U = xe^t[1 - x\{e^t - e^{-x(e^t-1)}\}]^{-1}.$$

We note that as $t \to \infty$, $U \to xe^t[1 - xe^t]^{-1}$ as would be expected since each $F_i(t) \to 1$. $U(x,t) = \sum_1^\infty F_n(t)(xe^t)^n$, so putting $y = xe^t$, we obtain

(5) $$F(y,t) = \sum_1^\infty F_n(t)y^n = y[1 - y + ye^{-t+y(e^{-t}-1)}]^{-1}.$$

Equation (5) implies

$$\sum_1^\infty (1 - F_n(t))y^n = \frac{y}{1-y}\frac{(e^{-t}y/(1-y))e^{(e^{-t}-1)y}}{1 + (e^{-t}y/(1-y))e^{(e^{-t}-1)y}}.$$

In order to find an asymptotic form for the mean time to stopping, $\mu_n$, we need to integrate this expression with respect to $t$. The term on the right-hand side is not directly integrable, but if we multiply by $1 + ye^{-t}$ it is. We have

$$\sum_1^\infty \int_0^\infty (1 + ye^{-t})(1 - F_n(t))y^n \, dt = \frac{y}{1-y}\left[-\ln\left\{1 + \frac{e^{-t}y}{1-y}e^{(e^{-t}-1)y}\right\}\right]_0^\infty$$

$$= -\frac{y}{1-y}\ln\{1-y\}$$

$$= (y + y^2 + y^3 + \cdots)\left(y + \frac{y^2}{2} + \frac{y^3}{3} + \cdots\right),$$

so that

$$\int_0^\infty (1 - F_n(t))\,dt + \int_0^\infty e^{-t}(1 - F_{n-1}(t))\,dt = 1 + \frac{1}{2} + \frac{1}{3} + \cdots + \frac{1}{n-1}.$$

Equation (8) shows that $F_{n-1}(t) \to 0$ as $n \to \infty$, so that the second integral above $\to 1$.

THEOREM 9. *The generating function for the duration of the annihilation process starting from $n$ particles in a block is*

$$F(y,t) = \sum_1^\infty F_n(t)y^n = y[1 - y + ye^{-t+y(e^{-t}-1)}]^{-1}.$$

*The mean time satisfies*

$$\left|\mu_n - \left(\frac{1}{2} + \frac{1}{3} + \cdots + \frac{1}{n-1}\right)\right| \to 0.$$

**Acknowledgment.** I am grateful to the referee for comments which have helped clarify several issues in this paper.

SCHOOL OF MATHEMATICAL SCIENCES
MONASH UNIVERSITY
AUSTRALIA 3800
E-MAIL: Aidan.Sudbury@sci.monash.edu.au